\documentclass[11pt,reqno]{amsart}
\usepackage{amsmath}
\usepackage{amssymb}
\usepackage{amstext}
\usepackage{amsfonts}
\usepackage{amsthm}
\usepackage{ifthen}
\usepackage{xspace}
\usepackage{graphicx}
\usepackage{mathrsfs}
\usepackage{enumerate}
\usepackage{breakcites}

\numberwithin{equation}{section}  

\newtheorem{punkt}{}[section]

\theoremstyle{plain}

\newtheorem{lemma}[punkt]{Lemma}
\newtheorem{proposition}[punkt]{Proposition}
\newtheorem{theorem}[punkt]{Theorem}

\theoremstyle{definition}

\theoremstyle{plain}
\newtheorem*{corollary*}{Corollary}
\newtheorem*{lemma*}{Lemma}
\newtheorem*{proposition*}{Proposition}
\newtheorem*{theorem*}{Theorem}

\theoremstyle{definition}
\newtheorem*{remark*}{Remark}
\newtheorem*{remarks*}{Remarks}
\newtheorem*{example*}{Example}
\newtheorem*{examples*}{Examples}
\newtheorem*{definition*}{Definition}
\newtheorem*{conjecture*}{Conjecture}
\newtheorem*{assumption*}{Assumption}
\newtheorem*{assumptions*}{Assumptions}
\newtheorem*{construction*}{Construction}

\def\mynat{\mathbb{N}}

\def\myreal{\mathbb{R}}

\def\myf{\mathcal{F}}

\def\ee{\mathbb{E}}
\def\pp{\mathbb{P}}

\def\re{\qopname\relax{no}{Re}\,}

\def\eg{e.g.\@\xspace}
\def\ie{i.e.\@\xspace}

\def\i{\qopname\relax{no}{\bf i}}

\begin{document}

\title{On the Second-Order Correlation Function \\[5pt] of the Characteristic Polynomial \\[5pt] of a Real Symmetric Wigner Matrix}

\author{H. K\"osters}
\address{Holger K\"osters, Fakult\"at f\"ur Mathematik, Universit\"at Bielefeld,
Postfach 100131, 33501 Bielefeld, Germany}
\email{hkoesters@math.uni-bielefeld.de}

\date{December 20, 2007}

\begin{abstract} 
We consider the asymptotic behaviour
of the second-order correlation function
of the characteristic polynomial 
of a real symmetric random matrix.
Our main result is that the existing result
for a random matrix from the Gaussian Orthogonal Ensemble,
obtained by \textsc{Br\'ezin} and \textsc{Hikami} \cite{BH2},
essentially continues to hold
for a general real symmetric Wigner matrix.
To~obtain this result, we adapt the approach
by \textsc{G\"otze} and \textsc{K\"osters} \cite{GK},
who proved the analogous result
for the Hermitian case.
\end{abstract}

\vspace*{0.5\baselineskip}

\maketitle

\markboth{H. K\"osters}{Characteristic Polynomials of Symmetric Wigner Matrices}

\section{Introduction}

In recent years, the characteristic polynomials 
of random matrices have found considerable interest.
This interest was sparked, at least in part, by the discovery 
by \textsc{Keating} and \textsc{Snaith} \cite{KS} that 
the moments of a random matrix from the Circular Unitary Ensemble (CUE)
seem to be related to
the moments of the Riemann zeta-function along the critical line.
Following this observation, several authors have investigated
the moments and correlation functions of the characteristic polynomial
also for other random matrix ensembles
(see \eg \textsc{Br\'ezin} and \textsc{Hikami} \cite{BH1,BH2},
\textsc{Mehta} and \textsc{Normand} \cite{MN},
\textsc{Strahov} and \textsc{Fyodorov} \cite{SF3},
\textsc{Baik}, \textsc{Deift} and \textsc{Strahov} \cite{BDS},
\textsc{Borodin} and \textsc{Strahov} \cite{BS},
\textsc{G\"otze} and \textsc{K\"osters} \cite{GK}).

One important observation is that the correlation functions
of the characteristic polynomials of Hermitian random matrices
are related to the ``sine kernel'' $\sin x / x$.
More precisely, this holds both for the unitary-invariant ensembles
(\textsc{Strahov} and \textsc{Fyodorov} \cite{SF3})
and --~at~least as far as the second-order correlation function is concerned~--
for the Hermitian Wigner ensembles
(\textsc{G\"otze} and \textsc{K\"osters} \cite{GK}).
Thus, the emergence of the sine kernel may be regarded
as ``universal'' for Hermitian random matrices.

In contrast to that, the correlation functions
of the characteristic polynomials of real symmetric random matrices
lead to different results.
This was first observed by \textsc{Br\'ezin} and \textsc{Hikami} \cite{BH2},
who investigated the Gaussian Orthogonal Ensemble (GOE)
(see \eg \textsc{Forrester} \cite{Fo} or \textsc{Mehta} \cite{Me} for definitions)
and came to the conclusion that the second-order correlation function
of the characteristic poly\-nomial is related to the function
$\sin x / x^3 - \cos x / x^2$ in this case.
(See below for a~more precise statement of this result.)
Moreover, \textsc{Borodin} and \textsc{Strahov} \cite{BS} 
obtained similar results for arbitrary products and ratios 
of characteristic poly\-nomials of the GOE.

\pagebreak[2]

The main purpose of this paper is to~generalize the above-mentioned result
by \textsc{Br\'ezin} and \textsc{Hikami} \cite{BH2}
about the second-order correlation function of the characteristic polynomial 
of the GOE to arbitrary real symmetric Wigner matrices.
Throughout this paper, we consider the following situation:
Let $Q$ be a prob\-ability distribution on the real line such that 
\begin{align}
\label{momentconditions} 
\int x \ Q(dx) = 0 \,,
\quad
a := \int x^2 \ Q(dx) = 1 \,,
\quad
b := \int x^4 \ Q(dx) < \infty \,,
\end{align}
and let $(X_{ii} / \sqrt{2})_{i \in \mynat}$
and $(X_{ij})_{i<j ,\, i,j \in \mynat}$
be independent families of independent real random variables
with distribution $Q$ on some probability space $(\Omega,\myf,\pp)$.
Also, let $X_{ji} := X_{ij}$ for $i < j$, $i,j \in \mynat$.
Then, for any $N \in \mynat$, 
the real symmetric Wigner matrix of size $N \times N$
is~given by $X_N = (X_{ij})_{1 \leq i,j \leq N}$,
and the second-order correlation function 
of the characteristic polynomial is~given by
\begin{align}
\label{correlationfunction}
f(N;\mu,\nu) := \ee(\det (X_N - \mu I_N) \cdot \det (X_N - \nu I_N)) \,,
\end{align}
where $\mu,\nu$ are real numbers 
and $I_N$ denotes the identity matrix of size $N \times N$.

In the special case where $Q$ is given by the standard Gaussian distribution,
the~distribution of the random matrix $X_N$ is
the Gaussian Orthogonal Ensemble (GOE).
(Note, however, that our scaling is slightly different
from that mostly used in the literature
(see \eg \textsc{Forrester} \cite{Fo} or \textsc{Mehta} \cite{Me}),
since the variance \linebreak[2] of the off-diagonal matrix entries
is fixed to $1$, and not to $1/2$.)
The result by \textsc{Br\'ezin} and \textsc{Hikami} \cite{BH2}
then corresponds to the statement that
\begin{multline}
\label{GOE}
\lim_{N \to \infty} \sqrt{\frac{2 \pi}{N^3}} \cdot \frac{1}{N!} \cdot e^{-N\xi^2/2} \cdot f \left( N; \sqrt{N} \xi + \frac{\mu}{\sqrt{N}\varrho(\xi)} , \sqrt{N} \xi + \frac{\nu}{\sqrt{N}\varrho(\xi)} \right) \qquad \\ = e^{\xi(\mu+\nu) / 2\varrho(\xi)} \cdot (2\pi\varrho(\xi))^3 \cdot \frac{1}{2} \left( \frac{\sin (\pi (\mu-\nu))}{(\pi (\mu-\nu))^3} - \frac{\cos (\pi (\mu-\nu))}{(\pi (\mu-\nu))^2} \right) \,,
\end{multline}
where $\xi \in (-2,+2)$, $\mu,\nu \in \mathbb{R}$,
and $\varrho(\xi) := \tfrac{1}{2\pi} \sqrt{4-\xi^2}$\,.

Our main result is the following generalization of (\ref{GOE}):

\begin{theorem}
\label{maintheorem}
Let $Q$ be a probability distribution on the real line
satisfying (\ref{momentconditions}),
let $f$ be defined as in (\ref{correlationfunction}),
let $\xi \in (-2,+2)$, and let $\mu,\nu \in \mathbb{R}$.
Then we~have
\begin{multline}
\label{symmetric}
\lim_{N \to \infty} \sqrt{\frac{2 \pi}{N^3}} \cdot \frac{1}{N!} \cdot e^{-N\xi^2/2} \cdot f \left( N; \sqrt{N} \xi + \frac{\mu}{\sqrt{N}\varrho(\xi)} , \sqrt{N} \xi + \frac{\nu}{\sqrt{N}\varrho(\xi)} \right) \qquad \\ = \exp \Big( \tfrac{b-3}{2} \Big) \cdot e^{\xi(\mu+\nu) / 2\varrho(\xi)} \cdot (2\pi\varrho(\xi))^3 \cdot T(\pi(\mu-\nu)) \,,
\end{multline}
where $\varrho(\xi) := \tfrac{1}{2\pi} \sqrt{4-\xi^2}$ 
denotes the density of the semi-circle law,
\begin{align}
\label{T}
T(x) := \frac{1}{2} \left( \frac{\sin x}{x^3} - \frac{\cos x}{x^2} \right) \quad\text{for}\quad x \ne 0 \,,
\end{align}
and $T(x) := 1/6$ for $x = 0$, by continuous extension.
\end{theorem}

In particular, we find that the correlation function 
of the characteristic poly\-nomial asymptotically factorizes
into a universal factor involving the function $T(x)$,
another universal factor involving the density $\varrho(\xi)$
of the semi-circle law, and \mbox{a~non-universal~factor}
which depends on the underlying distribution $Q$
only via its fourth moment $b$, or its fourth cumulant $b-3$.

It is interesting to compare the above result
with the corresponding result for Hermitian Wigner matrices 
(Theorem~1.1 in \textsc{G\"otze} and \textsc{K\"osters} \cite{GK}), 
which states that
under similar assumptions as in~(\ref{momentconditions})
and with similar notation as in~(\ref{correlationfunction}),
we have
\begin{multline}
\label{hermitian}
\lim_{N \to \infty} \sqrt{\frac{2\pi}{N}} \cdot \frac{1}{N!} \cdot e^{-N\xi^2/2} \cdot \hat{f} \left( N; \sqrt{N} \xi + \frac{\mu}{\sqrt{N}\varrho(\xi)} , \sqrt{N} \xi + \frac{\nu}{\sqrt{N}\varrho(\xi)} \right) \qquad \\ = \exp \Big( \hat{b} - \tfrac{3}{4} \Big) \cdot e^{\xi(\mu+\nu) / 2\varrho(\xi)} \cdot (2\pi\varrho(\xi)) \cdot \frac{\sin{\pi(\mu-\nu)}}{\pi(\mu-\nu)} \,.
\end{multline}
($\hat{f}$ and $\hat{b}$ denote the analogues of $f$ and $b$, respectively.)
Obviously, the structure of (\ref{hermitian}) and (\ref{symmetric}) is the same.
The most notable difference
is given by the fact that 
the ``sine kernel'' $\sin x / x$ in (\ref{hermitian})
is replaced with
the function $T(x)$ in (\ref{symmetric}).
It is noteworthy that both functions
are closely related to Bessel functions,
as~already observed by \textsc{Br\'ezin} and \textsc{Hikami} \cite{BH2}.
Indeed, it is well-known (see \eg p.\,78 in~\textsc{Erd\'elyi} \cite{Er2}) 
that
$$
\frac{\sin x}{x} = \sqrt{\pi} \, J_{1/2}(x) \,/\, (2x)^{1/2}
\quad\text{and}\quad
\frac{1}{2} \left( \frac{\sin x}{x^3} - \frac{\cos x}{x^2} \right) = \sqrt{\pi} \, J_{3/2}(x) \,/\, (2x)^{3/2} \,,
$$
where $J_p(x)$ denotes the Bessel function of order $p$.
Thus, if one wishes, one may rewrite
both (\ref{hermitian}) and (\ref{symmetric})
in the common form
\begin{multline*}
\lim_{N \to \infty} \frac{\sqrt{2\pi}}{N^p} \cdot \frac{1}{N!} \cdot e^{-N\xi^2/2} \cdot f \left( N; \sqrt{N} \xi + \frac{\mu}{\sqrt{N}\varrho(\xi)} , \sqrt{N} \xi + \frac{\nu}{\sqrt{N}\varrho(\xi)} \right) \qquad \\ = \sqrt{\pi} \cdot \exp \big( b^* \big) \cdot e^{\xi(\mu+\nu) / 2\varrho(\xi)} \cdot (2\pi\varrho(\xi))^{2p} \cdot \frac{J_p \big( \pi(\mu-\nu) \big)}{\big( 2\pi(\mu-\nu) \big)^p}
\end{multline*}
with $p := \tfrac{1}{2}$, $\hat{b}^* := \hat{b} - \tfrac{3}{4}$ in the Hermitian case
and $p := \tfrac{3}{2}$, $b^* := \frac{b-3}{2}$ in the symmetric case.

Furthermore, in the special case that $\xi = \mu = \nu = 0$,
Theorem \ref{maintheorem} reduces to a~result 
about determinants of random matrices, 
due to \textsc{Zhurbenko} \cite{Zh}.

To prove Theorem \ref{maintheorem}, we~show that 
the approach for Hermitian Wigner matrices
adopted by \textsc{G\"otze} and \textsc{K\"osters} \cite{GK}
can easily be adapted to real symmetric Wigner matrices.
This stands in contrast to the ``orthogonal polynomial approach''
typically used in the analysis of the invariant ensembles,
for which the transition 
from the unitary-invariant ensembles (such as the GUE)
to the orthogonal-invariant ensembles (such as the GOE)
is usually more complicated.

\medskip

\textbf{Acknowledgement.}
The author thanks Friedrich G\"otze for the suggestion
to~study the problem.

\bigskip

\section{Generating Functions}

In this section, we determine the exponential generating function
of the correlation function of the characteristic polynomial 
of a real symmetric Wigner matrix. 
Our results generalize those by \textsc{Zhurbenko} \cite{Zh},
who considered the special case of determinants.

We make the following conventions:
The determinant of the ``empty'' (\ie, $0 \times 0$) matrix 
is taken to be $1$. If $A$ is an $n \times n$ matrix 
and $z$ is a real or complex number, we set $A - z := A - z I_n$, 
where $I_n$ denotes the $n \times n$ identity matrix.
Also, if $A$ is an $n \times n$ matrix
and $i_1,\hdots,i_m$ and $j_1,\hdots,j_m$
are families of pairwise different indices
from the set $\{ 1,\hdots,n \}$,
we write $A^{[i_1,\hdots,i_m:j_1,\hdots,j_m]}$
for the $(n-m) \times (n-m)$-matrix obtained from $A$ 
by removing the rows indexed by $i_1,\hdots,i_m$
and the columns indexed by $j_1,\hdots,j_m$.
Thus, for any $n \times n$ matrix
$A = (a_{ij})_{1 \leq i,j \leq n}$ $(n \geq 1)$,
we~have
\begin{align}
\label{laplace-expansion}
\det(A) = \sum_{i,j=1}^{n-1} (-1)^{i+j-1} \, a_{i,n} \, a_{n,j} \, \det(A^{[n,i:n,j]}) + a_{n,n} \det(A^{[n:n]}) \,,
\end{align}
as follows by expanding the determinant
about the last row and the last column.
(For $n=1$, note that the big sum vanishes.)

Recall that we write $X_N$ for the real symmetric random matrix 
$(X_{ij})_{1 \leq i,j \leq N}$,
where the $X_{ij}$ are the random variables from the introduction.
We will analyze the~function
\begin{align*}
\ \ f(N;\mu,\nu) &:= \mathbb{E} \left( \det(X_N - \mu) \cdot \det(X_N - \nu) \right) \ \ \qquad\qquad & (N \geq 0) \,.
\end{align*}
To this purpose, we will also need the auxiliary functions
\begin{align*}
f_{11}^{A}(N;\mu,\nu) &:= \mathbb{E}(\det((X_N - \mu)^{[1:1]}) \cdot \det((X_N - \nu)^{[2:2]})) & (N \geq 2) \,,
\\
f_{11}^{B}(N;\mu,\nu) &:= \mathbb{E}(\det((X_N - \mu)^{[1:2]}) \cdot \det((X_N - \nu)^{[1:2]})) & (N \geq 2) \,,
\\
f_{11}^{C}(N;\mu,\nu) &:= \mathbb{E}(\det((X_N - \mu)^{[1:2]}) \cdot \det((X_N - \nu)^{[2:1]})) & (N \geq 2) \,,
\allowdisplaybreaks\\
f_{10}(N;\mu,\nu) &:= \mathbb{E}(\det(X_{N-1} - \mu) \cdot \det(X_N - \nu)) & (N \geq 1) \,,
\\
f_{01}(N;\mu,\nu) &:= \mathbb{E}(\det(X_N - \mu) \cdot \det(X_{N-1} - \nu)) & (N \geq 1) \,.
\end{align*}
Note that the functions $f_{11}^B$ and $f_{11}^C$ actually coincide,
but we will not need this.
Since $\mu$ and $\nu$ can be regarded as constants
for the purposes of this section, we~will only write
$f(N)$ instead of $f(N;\mu,\nu)$, etc.

\pagebreak[2]

We have the following recursive equations:

\begin{lemma}
\label{fn-recursion}
\begin{align}
f(0) &= 1 \,, \nonumber \\
f(N) &= (2 + \mu \nu) \, f(N-1) + b (N-1) \, f(N-2)
\nonumber\\&\qquad \,+\, (N-1)(N-2) \, f_{11}^A(N-1)
\nonumber\\&\qquad \,+\, (N-1)(N-2) \, f_{11}^B(N-1)
\nonumber\\&\qquad \,+\, (N-1)(N-2) \, f_{11}^C(N-1)
\nonumber\\&\qquad \,+\, \nu(N-1) \, f_{10}(N-1)
\nonumber\\&\qquad \,+\, \mu(N-1) \, f_{01}(N-1) & (N \geq 1) \,,
\label{fn-recursion-1}
\allowdisplaybreaks\\[+5pt]
f_{11}^A(N) &= \mu \nu \, f(N-2) + (N-2) \, f(N-3)
\nonumber\\&\qquad \,+\, (N-2)(N-3) \, f_{11}^A(N-2)
\nonumber\\&\qquad \,+\, \nu(N-2) \, f_{10}(N-2)
\nonumber\\&\qquad \,+\, \mu(N-2) \, f_{01}(N-2) & (N \geq 2) \,,
\label{fn-recursion-2}
\allowdisplaybreaks\\[+5pt]
f_{11}^B(N) &= f(N-2) + (N-2) \, f(N-3)
\nonumber\\&\qquad \,+\, (N-2)(N-3) \, f_{11}^B(N-2) & (N \geq 2) \,,
\label{fn-recursion-3}
\allowdisplaybreaks\\[+5pt]
f_{11}^C(N) &= f(N-2) + (N-2) f(N-3)
\nonumber\\&\qquad \,+\, (N-2)(N-3) \, f_{11}^C(N-2) & (N \geq 2) \,,
\label{fn-recursion-4}
\allowdisplaybreaks\\[+5pt]
f_{10}(N) &= -(N-1) \, f_{01}(N-1) - \nu \, f(N-1) & (N \geq 1) \,,
\label{fn-recursion-5}
\allowdisplaybreaks\\[+5pt]
f_{01}(N) &= -(N-1) \, f_{10}(N-1) - \mu \, f(N-1) & (N \geq 1) \,.
\label{fn-recursion-6}
\end{align}
\end{lemma}

\begin{proof}
We give the proof for the recursive equation for $f(N)$ only,
the proofs for the remaining recursive equations being very similar.

The result for $f(0)$ is clear.
For $N \geq 1$, we expand the determinants 
of the matrices $(X_N-\mu)$ and $(X_N-\nu)$
as in (\ref{laplace-expansion}) and use the independence
of the random variables $X_{ij} = X_{ji}$ ($i \leq j$),
thereby obtaining
\begin{align*}
&\mskip24mu f(N) \\
&= \sum_{i,j=1}^{N-1} \sum_{k,l=1}^{N-1} (-1)^{i+j+k+l} \, \ee \left( X_{i,N} X_{N,j} X_{k,N} X_{N,l} \right) \cdot \ee \left( \det(X_{N-1}-\mu)^{[i:j]} \cdot \det(X_{N-1}-\nu)^{[k:l]} \right) \\
     &\quad \,+\, \sum_{i,j=1}^{N-1} (-1)^{i+j+1} \, \ee \left( X_{i,N} X_{N,j} \right) \cdot \ee \left( X_{N,N}-\nu \right) \cdot \ee \left(\det(X_{N-1}-\mu)^{[i:j]} \cdot \det(X_{N-1}-\nu) \right) \\
     &\quad \,+\, \sum_{k,l=1}^{N-1} (-1)^{k+l+1} \, \ee \left( X_{k,N} X_{N,l} \right) \cdot \ee \left( X_{N,N}-\mu \right) \cdot \ee \left(\det(X_{N-1}-\mu) \cdot \det(X_{N-1}-\nu)^{[k:l]} \right) \\
     &\quad \,+\, \ee \left( (X_{N,N}-\mu) (X_{N,N}-\nu) \right) \cdot \ee \left( \det(X_{N-1}-\mu) \cdot \det(X_{N-1}-\nu) \right) \,.
\end{align*}
Since the random variables $X_{ij} = X_{ji}$ ($i \le j$) are independent 
with $\ee(X_{ij}) = 0$ ($i \leq j$), several of the expectations vanish, 
and the sum reduces to
\begin{align*}
&\mskip24mu f(N) \\
&= ( \ee X_{N,N}^2 + \mu \nu ) \cdot \ee \left( \det(X_{N-1} - \mu) \cdot \det(X_{N-1} - \nu) \right) \\
     &\quad \,+\, \sum_{i = j  =  k = l} \ee X_{i,N}^4 \cdot \ee \left( \det(X_{N-1}-\mu)^{[i:j]} \cdot \det(X_{N-1}-\nu)^{[k:l]} \right) \\
     &\quad \,+\, \sum_{i = j \ne k = l} \ee X_{i,N}^2 \cdot \ee X_{k,N}^2 \cdot \ee \left( \det(X_{N-1}-\mu)^{[i:j]} \cdot \det(X_{N-1}-\nu)^{[k:l]} \right) \\
     &\quad \,+\, \sum_{i = k \ne j = l} \ee X_{i,N}^2 \cdot \ee X_{j,N}^2 \cdot \ee \left( \det(X_{N-1}-\mu)^{[i:j]} \cdot \det(X_{N-1}-\nu)^{[k:l]} \right) \\
     &\quad \,+\, \sum_{i = l \ne j = k} \ee X_{i,N}^2 \cdot \ee X_{k,N}^2 \cdot \ee \left( \det(X_{N-1}-\mu)^{[i:j]} \cdot \det(X_{N-1}-\nu)^{[k:l]} \right) \\
     &\quad \,+\, \nu \sum_{i=j} \ee X_{i,N}^2 \cdot \ee \left(\det(X_{N-1}-\mu)^{[i:j]} \cdot \det(X_{N-1}-\nu) \right) \\
     &\quad \,+\, \mu \sum_{k=l} \ee X_{k,N}^2 \cdot \ee \left(\det(X_{N-1}-\mu) \cdot \det(X_{N-1}-\nu)^{[k:l]} \right) \,.
\end{align*}
From this (\ref{fn-recursion-1}) follows by noting that
$\ee X_{N,N}^2 = 2$,
$\ee X_{i,N}^2 = 1$,
$\ee X_{i,N}^4 = b$,
and by exploiting obvious symmetries.
\end{proof}

It turns out that the above recursions can be combined 
into a single recursion involving only the values $f(N)$.
Using the abbreviations
$$
c(N) := \frac{f(N)}{N!} \qquad (N \geq 0)
$$
and \mbox{\qquad\qquad\qquad}
$$
s(N) := \sum_{\substack{k=0,\hdots,N \\ \text{$k$ even}}} c(N-k) \qquad (N \geq 0) \,,
$$
we have the following result:

\begin{lemma}
\label{fn-superrecursion}
The values $c(N)$ satisfy the recursive equation
\begin{align}
\label{cn-formula-3}
c(0) &= 1 \,, \\
N c(N) &= 2 \cdot c(N-1) + (N+1) \cdot c(N-2)
    \nonumber\\&\qquad \,+\, \mu \nu \cdot \big( s(N-1) + s(N-3) \big)
    \nonumber\\&\qquad \,-\, (\mu^2 + \nu^2) \cdot s(N-2)
    \nonumber\\&\qquad \,+\, (b-3) \cdot \big( c(N-2) - c(N-4) \big) & (N \geq 1) \,,
\label{cn-formula-4}
\end{align}
where all terms $c(\,\cdot\,)$ and $s(\,\cdot\,)$
with a negative argument are taken to be zero.
\end{lemma}

\begin{proof}
It follows from Lemma \ref{fn-recursion} that
$$
f(N-2) = f_{11}^A(N-1) + f_{11}^B(N-1) + f_{11}^C(N-1) + (b-3) (N-3) \, f(N-4)
$$
for all $N \geq 3$. 
Thus, we~can substitute $f_{11}^A(N-1) + f_{11}^B(N-1) + f_{11}^C(N-1)$
on the right-hand side of (\ref{fn-recursion-1}) to~obtain
\begin{align*}
f(N) &= (2 + \mu \nu) \, f(N-1) + b(N-1) \, f(N-2)
    \\&\quad \,+\, (N-1)(N-2) \, \Big( f(N-2) - (b-3) (N-3) f(N-4) \Big)
    \\&\quad \,+\, \nu(N-1) \, f_{10}(N-1)
    \\&\quad \,+\, \mu(N-1) \, f_{01}(N-1)
\\
&= (2 + \mu \nu) \, f(N-1) + (N+1) (N-1) \, f(N-2)
     \\&\quad \,+\, \nu(N-1) \, f_{10}(N-1)
     \\&\quad \,+\, \mu(N-1) \, f_{01}(N-1) 
     \\&\quad \,+\, (b-3) \cdot \Big( (N-1) \, f(N-2) - (N-1)(N-2)(N-3) \, f(N-4) \Big)
\end{align*}
for all $N \geq 3$. Dividing by $(N-1)!$, it follows that
\begin{align*}
N c(N) &= (2 + \mu \nu) \cdot c(N-1) + (N+1) \, c(N-2)
     \\&\quad \,+\, \nu \cdot f_{10}(N-1) \,/\, (N-2)!
     \\&\quad \,+\, \mu \cdot f_{01}(N-1) \,/\, (N-2)!
     \\&\quad \,+\, (b-3) \cdot \Big( c(N-2) - c(N-4) \Big)
\end{align*}
for all $N \geq 3$. (For $N=3$, note that the second term 
in the large bracket \pagebreak[1] vanishes.)
Since
\begin{align*}
f_{10}(N-1) \,/\, (N-2)! &= - \nu s(N-2) + \mu s(N-3) \,,
\\[+5pt]
f_{01}(N-1) \,/\, (N-2)! &= - \mu s(N-2) + \nu s(N-3) \,,
\end{align*}
for all $N \geq 3$, as follows from (\ref{fn-recursion-5}) and (\ref{fn-recursion-6})
by a straightforward induction, the~assertion for $N \geq 3$ is proved.

\pagebreak[1]

The~assertion for $N < 3$ follows from Lemma \ref{fn-recursion}
by direct calculation:
\begin{align*}
  c(0) = f(0) &= 1 
\allowdisplaybreaks\\[+5pt]
1 c(1) = f(1) &= (2 + \mu \nu) f(0)
               = (2 + \mu \nu) c(0)
               = 2 c(0) + \mu \nu s(0) 
\allowdisplaybreaks\\[+5pt]
2 c(2) = f(2) &= (2 + \mu\nu) f(1) + b f(0) + \nu (-\nu f(0)) + \mu (-\mu f(0))
\\             &= (2 + \mu\nu) f(1) + b f(0) - (\mu^2 + \nu^2) f(0)
\\             &= (2 + \mu\nu) c(1) + b \mskip1mu c(0) - (\mu^2 + \nu^2) c(0)
\\             &= 2 c(1) + 3 c(0) + \mu\nu c(1) - (\mu^2 + \nu^2) c(0) + (b-3) c(0)
\\             &= 2 c(1) + 3 c(0) + \mu\nu s(1) - (\mu^2 + \nu^2) s(0) + (b-3) c(0) 
\end{align*}
\end{proof}

Using Lemma \ref{fn-superrecursion}, we can determine
the exponential generating function of the~sequence 
$(f(N))_{N \geq 0}$:

\begin{lemma}
\label{fn-genfun}
The exponential generating function $F(x) := \sum_{N=0}^{\infty} f(N) \, x^N \,/\, N!$ 
of the sequence $(f(N))_{N \geq 0}$ is given by
$$
F(x) = \frac{\exp \left( \mu \nu \cdot \frac{x}{1-x^2} - \tfrac{1}{2} (\mu^2 + \nu^2) \cdot \frac{x^2}{1-x^2} + b^* x^2 \right)}{(1-x)^{5/2} \cdot (1+x)^{1/2}} \,,
$$
where $b^* := \tfrac{1}{2}(b - 3)$.
\end{lemma}

\begin{proof}
Multiplying (\ref{cn-formula-4}) by $x^{N-1}$, summing over $N$
and recalling our convention con\-cerning negative arguments, we have
\begin{align*}
  \sum_{N=1}^{\infty} N c(N) x^{N-1}
= \sum_{N=1}^{\infty} & 2 c(N-1) x^{N-1}
+ \sum_{N=2}^{\infty} (N+1) c(N-2) x^{N-1}
\\&+ \mu \nu \left( \sum_{N=1}^{\infty} s(N-1) x^{N-1} + \sum_{N=3}^{\infty} s(N\!-\!3) x^{N-1} \right)
\\&- (\mu^2+\nu^2) \sum_{N=2}^{\infty} s(N-2) x^{N-1}
\\&+ 2b^* \left( \sum_{N=2}^{\infty} c(N-2) x^{N-1} - \sum_{N=4}^{\infty} c(N-4) x^{N-1} \right) \,,
\end{align*}
whence
\begin{multline*}
F'(x) = 2 F(x) + (3x F(x) + x^2 F'(x)) \\ \,+\, \mu \nu \frac{1+x^2}{1-x^2} F(x) - (\mu^2+\nu^2) \frac{x}{1-x^2} F(x) + 2b^* \left( x F(x) - x^3 F(x) \right) \,.
\end{multline*}
This leads to the differential equation
$$
F'(x) = \left( \frac{2 + 3x}{1-x^2} + \mu \nu \frac{1+x^2}{(1-x^2)^2} - (\mu^2+\nu^2) \frac{x}{(1-x^2)^2} + 2b^*x \right) F(x) \,,
$$
which has the solution
$$
F(x) = \frac{F_0}{(1-x)^{5/2} \cdot (1+x)^{1/2}} \exp \left( \mu \nu \frac{x}{1-x^2} - \tfrac{1}{2} (\mu^2 + \nu^2) \frac{1}{1-x^2} + b^* x^2 \right) \,.
$$
Here, $F_0$ is a multiplicative constant
which must be chosen as $\exp(\tfrac{1}{2} (\mu^2 + \nu^2))$
in~order to satisfy (\ref{cn-formula-3}).
Inserting this above completes the proof.
\end{proof}

\bigskip

\section{The Proof of Theorem \ref{maintheorem}}

This section is devoted to proving Theorem \ref{maintheorem}.
In doing so, we will closely follow the proof of Theorem~1.1
in \textsc{G\"otze} and \textsc{K\"osters} \cite{GK}.
Throughout this~section, $T(x)$ will denote the function 
defined in Theorem \ref{maintheorem}.

We will first establish the following slightly more general result:

\begin{proposition}
\label{proposition-B}
Let $Q$ be a probability distribution on the real line
satisfying~(\ref{momentconditions}),
let $f$ be defined as in (\ref{correlationfunction}),
let $(\xi_N)_{N \in \mathbb{N}}$ be a sequence of real numbers
such~that $\lim_{N \to \infty} \xi_N / \sqrt{N} = \xi$
for some $\xi \in (-2,+2)$,
and let $\eta \in \mathbb{C}$.
Then we~have
\begin{multline*}
\lim_{N \to \infty} \sqrt{\frac{2\pi}{N^3}} \cdot \frac{1}{N!} \cdot \exp(-\xi_N^2/2) \cdot f \left( N;\xi_N+\frac{\eta}{\sqrt{N}},\xi_N-\frac{\eta}{\sqrt{N}} \right)
\\
= \exp \left( \tfrac{b-3}{2} \right) \cdot (4-\xi^2)^{3/2} \cdot T \big( \sqrt{4-\xi^2} \cdot \eta \big) \,.
\end{multline*}
\end{proposition}

\pagebreak[2]

It is easy to see that Proposition~\ref{proposition-B} implies Theorem~\ref{maintheorem}:

\begin{proof}[Proof of Theorem \ref{maintheorem}]
Taking
$$
\xi_N := \sqrt{N} \xi + \frac{\pi(\mu+\nu)}{\sqrt{N} \cdot \sqrt{4-\xi^2}}
\qquad\text{and}\qquad
\eta := \frac{\pi(\mu-\nu)}{\sqrt{4-\xi^2}}
$$
in Proposition \ref{proposition-B}, we have
\begin{multline*}
\lim_{N \to \infty} \sqrt{\frac{2\pi}{N^3}} \cdot \frac{1}{N!} \cdot \exp \left(- N \xi^2 / 2 - 2 \pi \xi (\mu+\nu) / 2 \sqrt{4-\xi^2} \right) \\ \,\cdot\, f \left( N;\xi_N+\frac{2\pi\mu}{\sqrt{N} \sqrt{4-\xi^2}},\xi_N+\frac{2\pi\nu}{\sqrt{N} \sqrt{4-\xi^2}} \right)
\\
= \exp \left( \tfrac{b-3}{2} \right) \cdot (4-\xi^2)^{3/2} \cdot T \big( \pi(\mu-\nu) \big) \,,
\end{multline*}
from which Theorem \ref{maintheorem} follows by a simple rearrangement.
\end{proof}

\begin{proof}[Proof of Proposition \ref{proposition-B}]
From the exponential generating function obtained in Lemma \ref{fn-genfun}, 
we have the integral representation
\begin{align}
\label{intrep1}
\frac{f(N;\mu,\nu)}{N!}
=
\frac{1}{2 \pi i} \int_\gamma \frac{\exp \left( \mu \nu \cdot \frac{z}{1-z^2} - \tfrac{1}{2} (\mu^2 + \nu^2) \cdot \frac{z^2}{1-z^2} + b^* z^2 \right)}{(1-z)^{5/2} \cdot (1+z)^{1/2}} \ \frac{dz}{z^{N+1}} \,,
\end{align}
where $\gamma \equiv \gamma_N$ denotes the counterclockwise circle of radius 
$R \equiv R_N = 1 - 1/N$ around the origin.
(We will assume that $N \geq 2$ throughout the~proof.) 

Setting $\mu = \xi_N + \eta / \sqrt{N}$ and $\nu = \xi_N - \eta / \sqrt{N}$
and doing a simple calculation, we~obtain
\begin{align*}
&\mskip24mu 
   \exp \left( \mu \nu \cdot \frac{z}{1-z^2} - \tfrac{1}{2} (\mu^2 + \nu^2) \cdot \frac{z^2}{1-z^2} + b^* z^2 \right) \\
&= \exp \left( (\xi_N^2 - \eta^2 / N) \cdot \frac{z}{1-z^2} - (\xi_N^2 + \eta^2 / N) \cdot \frac{z^2}{1-z^2} + b^* z^2  \right) \\
&= \exp \left( \tfrac{1}{2} \xi_N^2 + \eta^2 / N \right) \cdot \exp \left( - \tfrac{1}{2} \xi_N^2 \cdot \frac{1-z}{1+z} - (\eta^2 / N) \cdot \frac{1}{1-z} + b^* z^2  \right) \,.
\end{align*}
Thus,
\begin{align}
\label{intrep2}
\frac{1}{N!} \cdot f & \left( N;\xi_N+\frac{\eta}{\sqrt{N}},\xi_N-\frac{\eta}{\sqrt{N}} \right)
=
\exp \left( \tfrac{1}{2} \xi_N^2 + \eta^2 / N \right)
\nonumber\\& \,\cdot\, 
\frac{1}{2 \pi i} \int_\gamma \frac{\exp \left( - \tfrac{1}{2} \xi_N^2 \cdot \frac{1-z}{1+z} - (\eta^2 / N) \cdot \frac{1}{1-z} + b^* z^2  \right)}{(1-z)^{5/2} \cdot (1+z)^{1/2}} \ \frac{dz}{z^{N+1}} \,.
\end{align}
Putting
$$
h(z) := \frac{\exp \left( - \tfrac{1}{2} \xi_N^2 \cdot \frac{1-z}{1+z} - (\eta^2 / N) \cdot \frac{1}{1-z} + b^* z^2  \right)}{(1-z)^{5/2} \cdot (1+z)^{1/2}}
$$
and
$$
h_0(z) := \frac{\exp(b^*)}{\sqrt{2}} \cdot \frac{\exp \left( - \tfrac{1}{4} \xi_N^2 \cdot (1-z) - (\eta^2 / N) \cdot \frac{1}{1-z} \right)}{(1-z)^{5/2}}
$$
we can rewrite the integral in (\ref{intrep2}) as
$$
\frac{1}{2 \pi i} \int_\gamma h(z) \ \frac{dz}{z^{N+1}}
= I_1 + I_2 + I_3 - I_4 \,,
$$
where
\begin{align*}
I_1 &:= \frac{1}{2 \pi i} \int_\gamma h_0(z) \ \frac{dz}{z^{N+1}} \,,
\\
I_2 &:= \frac{1}{2 \pi} \int_{1/\sqrt{N}}^{2\pi-1/\sqrt{N}} h(Re^{it}) \ \frac{dt}{(Re^{it})^{N}} \,,
\\
I_3 &:= \frac{1}{2 \pi} \int_{-1/\sqrt{N}}^{+1/\sqrt{N}} \left( h(Re^{it}) - h_0(Re^{it}) \right) \ \frac{dt}{(Re^{it})^{N}} \,,
\\
I_4 &:= \frac{1}{2 \pi} \int_{1/\sqrt{N}}^{2\pi-1/\sqrt{N}} h_0(Re^{it}) \ \frac{dt}{(Re^{it})^{N}} \,.
\end{align*}
We will show that the integral $I_1$ is the asymptotically dominant term.

First~of~all, note that since $\xi_N \in \myreal$, we have
\begin{align}
\label{xi-est}
\left| \exp \big( - \tfrac{1}{4} \xi_N^2 \cdot (1-z) \big) \right| = \exp \big( - \tfrac{1}{4} \xi_N^2 \cdot \re(1-z) \big) \leq 1
\end{align}
for any $z \in \mathbb{C}$ with $\re(z) \leq 1$.
Using the series expansion
$$
  \exp \left( - (\eta^2 / N) \cdot \frac{1}{1-z} \right)
= \sum_{l=0}^{\infty} \frac{(-1)^l \, \eta^{2l}}{l! \, N^l} \, \frac{1}{(1-z)^l} \,,
$$
which converges uniformly on the contour $\gamma$
(for fixed $N \geq 2$ and fixed $\eta \in \mathbb{C}$), 
we~obtain
\begin{align}
\label{I1}
\frac{I_1}{N^{3/2}} = \frac{\exp(b^*)}{\sqrt{2}} \cdot \sum_{l=0}^{\infty} \frac{(-1)^l \, \eta^{2l}}{l! \, N^{l+3/2}} \cdot \frac{1}{2 \pi i} \int_{\gamma} \frac{\exp \left( - \tfrac{1}{4} \xi_N^2 \cdot (1-z) \right)}{(1-z)^{l+5/2}} \ \frac{dz}{z^{N+1}} \,.
\end{align}

For each $l=0,1,2,3,\hdots,$ the integral in (\ref{I1}) may be rewritten as
\begin{multline}
\frac{1}{2 \pi i} \int_{\gamma} \frac{\exp \left( - \tfrac{1}{4} \xi_N^2 \cdot (1-z) \right)}{(1-z)^{l+5/2}} \ \frac{dz}{z^{N+1}}
\\
= \frac{1}{2 \pi i} \int_{R-\i\infty}^{R+\i\infty} \frac{\exp \left( - \tfrac{1}{4} \xi_N^2 \cdot (1-z) \right)}{(1-z)^{l+5/2}} \ \frac{dz}{z^{N+1}} \,.
\label{I2}
\end{multline}
Indeed, for any $R' > 1$, we can replace
the contour $\gamma$ by the contour $\delta$
consisting of the line segment between the points
$R - \i \sqrt{(R')^2 - R^2}$ and $R + \i \sqrt{(R')^2 - R^2}$,
and the arc of radius $R'$ around the origin to the left of this line segment.
%
%
It is easy to see then that the integral along the~arc is~bounded by
$$
\frac{1}{2\pi} \cdot 2\pi R' \cdot \frac{1}{(R'-1)^{l+5/2}} \cdot \frac{1}{(R')^{N+1}} \,,
$$
which tends to zero as $R' \to \infty$. This proves (\ref{I2}).
Next, performing a change of~variables, we obtain, for each $l=0,1,2,3,\hdots,$
\begin{multline}
\frac{1}{2 \pi i} \int_{R-\i\infty}^{R+\i\infty} \frac{\exp \left( - \tfrac{1}{4} \xi_N^2 \cdot (1-z) \right)}{(1-z)^{l+5/2}} \ \frac{dz}{z^{N+1}} \\
= N^{l+3/2} \cdot \frac{1}{2 \pi} \int_{-\infty}^{+\infty} \frac{\exp \left( - \tfrac{1}{4} (\xi_N^2/N) \cdot (1-iu) \right)}{(1-iu)^{l+5/2}} \ \frac{du}{(1-\frac{1-iu}{N})^{N+1}} \,.
\label{I3}
\end{multline}
Since $\lim_{N \to \infty} \xi_N / \sqrt{N} = \xi$,
an application of the dominated convergence theorem yields, 
for each $l=0,1,2,3,\hdots,$
\begin{multline}
\lim_{N \to \infty} \frac{1}{2 \pi} \int_{-\infty}^{+\infty} \frac{\exp \left( - \tfrac{1}{4} (\xi_N^2/N) \cdot (1-iu) \right)}{(1-iu)^{l+5/2}} \ \frac{du}{(1-\frac{1-iu}{N})^{N+1}} \\ 
= \frac{1}{2 \pi} \int_{-\infty}^{+\infty} \frac{\exp \left( (1 - \tfrac{1}{4} \xi^2) \cdot (1-iu) \right)}{(1-iu)^{l+5/2}} \ du \,.
\label{I4}
\end{multline}
Using the Laplace inversion formula (see \eg Chapter~24 in \textsc{Doetsch} \cite{Do}) 
and the~functional equation of the~Gamma function, it follows that the final integral
is~equal to
$$
  \frac{(1-\tfrac{1}{4}\xi^2)^{l+3/2}}{\Gamma(l+5/2)}
= \frac{1}{\sqrt{\pi}} \cdot \frac{(l+1)!}{(2l+3)!} \cdot (4-\xi^2)^{l+3/2} \,.
$$
Thus, putting it all together, we have shown that for each $l=0,1,2,3,\hdots,$
\begin{multline}
\label{term-by-term}
\lim_{N \to \infty} \left( \frac{(-1)^l \,\eta^{2l}}{l! \, N^{l+3/2}} \cdot \frac{1}{2 \pi i} \int_{\gamma} \frac{\exp \left( - \tfrac{1}{4} \xi_N^2 \cdot (1-z) \right)}{(1-z)^{l+5/2}} \ \frac{dz}{z^{N+1}} \right) \\
= \frac{1}{\sqrt{\pi}} \cdot (4-\xi^2)^{3/2} \cdot \frac{(-1)^l \, (l+1) \, ( \sqrt{4 - \xi^2} \cdot \eta)^{2l}}{(2l+3)!} \,.
\end{multline}
In particular, the series in (\ref{I1}) converges termwise.

\pagebreak[2]

We will show that the series in (\ref{I3}) also converges as a whole.
To this purpose, let $\varepsilon^2 > 0$ denote a positive constant 
such that $\cos t \leq 1 - \varepsilon^2 t^2$ for $-\pi \leq t \leq +\pi$.
Then, for any $\alpha > 0$ and any $-\pi \leq t_1 < t_2 \leq +\pi$, 
we have the estimate
\begin{align}
      \int_{t_1}^{t_2} \frac{1}{|1-Re^{it}|^{\alpha}} \ dt
&=    \int_{t_1}^{t_2} \frac{1}{(1+R^2-2R \cos t)^{\alpha/2}} \ dt \nonumber\\
&\leq \int_{t_1}^{t_2} \frac{1}{((1-R)^2+\varepsilon^2 t^2)^{\alpha/2}} \ dt \nonumber\\
&=    N^\alpha \int_{t_1}^{t_2} \frac{1}{(1+N^2 \varepsilon^2 t^2)^{\alpha/2}} \ dt \nonumber\\
&=    KN^{\alpha-1} \int_{N\varepsilon t_1}^{N\varepsilon t_2} \frac{1}{(1+u^2)^{\alpha/2}} \ du \,,
\label{theestimate}
\end{align}
where $K$ denotes some absolute constant.
Let us convene that this constant $K$ may change 
from occurrence to occurrence in the subsequent calculations. 
Then it~follows from (\ref{xi-est}) and (\ref{theestimate})
that, for each $l=0,1,2,3,\hdots,$
\begin{align*}
\Bigg| \frac{(-1)^l \, \eta^{2l}}{l! \, N^{l+3/2}} & \cdot \frac{1}{2 \pi i} \int_{\gamma} \frac{\exp \left( - \tfrac{1}{4} \xi^2_N \cdot (1-z) \right)}{(1-z)^{l+5/2}} \ \frac{dz}{z^{N+1}} \Bigg|
\\&\leq \frac{|\eta|^{2l}}{l!} \cdot \frac{1}{N^{l+3/2}} \cdot \frac{1}{2\pi} \int_{-\pi}^{+\pi} \frac{1}{|1-Re^{it}|^{l+5/2}} \ \frac{dt}{|Re^{it}|^N} 
\\&\leq K \cdot \frac{|\eta|^{2l}}{l!} \cdot \left( 1 + \int_{1}^{\infty} \frac{1}{u^{l+5/2}} \ du \right)
\\&\leq K \cdot \frac{|\eta|^{2l}}{l!} \cdot \left( 1 + \frac{1}{l+3/2} \right) \,,
\end{align*}
where
\begin{align*}
K \cdot \sum_{l=0}^{\infty} \frac{|\eta|^{2l}}{l!} \cdot \left( 1 + \frac{1}{l+3/2} \right) < \infty \,.
\end{align*}
It therefore follows from (\ref{term-by-term}) that
\begin{align*}
   \lim_{N\to\infty} \frac{I_1}{N^{3/2}} 
&= \frac{\exp(b^*)}{\sqrt{2\pi}} \cdot (4-\xi^2)^{3/2} \cdot \sum_{l=0}^{\infty} \frac{(-1)^l (l+1) \big( \sqrt{4-\xi^2} \cdot \eta \big)^{2l}}{(2l+3)!} \\
&= \frac{\exp(b^*)}{\sqrt{2\pi}} \cdot (4-\xi^2)^{3/2} \cdot T(\sqrt{4-\xi^2} \cdot \eta) \,,
\end{align*}
since
$$
\sum_{l=0}^{\infty} \frac{(-1)^l \, (l+1) \, z^{2l}}{(2l+3)!} = \frac{1}{2} \left( \frac{\sin z}{z^3} - \frac{\cos z}{z^2} \right) = T(z) \,.
$$
Hence, to complete the proof of Proposition \ref{proposition-B}, it remains 
to show that \linebreak $\lim_{N \to \infty} I_j / N^{3/2} = 0$; $j=2,3,4$.

For the integral $I_2$, we use the estimates
($R = 1 - 1/N$, $t \in \mathbb{R}$)
\begin{align}
   \left| \exp \left( - \tfrac{1}{2} \xi_N^2 \cdot \frac{1-Re^{it}}{1+Re^{it}} \right) \right|
&= \exp \left( - \tfrac{1}{2} \xi_N^2 \cdot \re \left( \frac{1-Re^{it}}{1+Re^{it}} \right) \right) \nonumber\\
&= \exp \left( - \tfrac{1}{2} \xi_N^2 \cdot \frac{1-R^2}{1+R^2+2R\cos t} \right)
\leq 1 \,,
\label{est4}
\\[+7pt]
      \left| \exp \left( - (\eta^2 / N) \cdot \frac{1}{1-Re^{it}} \right) \right|
&\leq \exp \left( (|\eta|^2 / N) \cdot \frac{1}{|1-Re^{it}|} \right) \nonumber\\
&\leq \exp \left( (|\eta|^2 / N) \cdot \frac{1}{1-R} \right)
   =  \exp(|\eta|^2) \,,
\label{est5}
\\[+7pt]
\label{est6}
      \left| \exp \Big( b^* (Re^{it})^2 \Big) \right|
&\leq \exp \Big( |b^*| |Re^{it}|^2 \Big) 
 \leq \exp (|b^*|)
\end{align}
as well as (\ref{theestimate}) to obtain
\begin{align*}
|I_2|
&\leq \frac{1}{2\pi} \int_{1/\sqrt{N}}^{2\pi-1/\sqrt{N}} \frac{\exp(|\eta|^2 + |b^*|)}{|1-Re^{it}|^{5/2} \cdot |1+Re^{it}|^{1/2}} \ \frac{dt}{R^N} \\
&\leq K \exp(|\eta|^2 + |b^*|) \cdot 
\left(
  \int_{1/\sqrt{N}}^{\pi/2} \frac{1}{|1-Re^{it}|^{5/2}} \ dt
+ \int_{0}^{\pi/2} \frac{1}{|1-Re^{it}|^{1/2}} \ dt
\right) \\
&\leq K \exp(|\eta|^2 + |b^*|) \cdot
\left(
  N^{3/2} \int_{\sqrt{N}\varepsilon}^{\infty} \frac{1}{u^{5/2}} \ du
+ N^{-1/2} \left( 1 + \int_{1}^{N\varepsilon\pi/2} \frac{1}{u^{1/2}} \ du \right)
\right) \\
&\leq K \exp(|\eta|^2 + |b^*|) \cdot
\left(
  N^{3/4} 
+ 1
\right)
= o(N^{3/2}) \,,
\end{align*}
where $K$ denotes some absolute constant which may change from step to step. 

For the integral $I_3$, we write
$$
h(z) - h_0(z) =
\frac{\exp \left( - \tfrac{1}{4} \xi_N^2 \cdot (1-z) - (\eta^2 / N) \cdot \frac{1}{1-z} \right)}{(1-z)^{5/2}}
\cdot
\left( \tilde{h}(z) - \tilde{h}(1) \right) \,,
$$
where
$$
\tilde{h}(z) = \frac{\exp \left( - \tfrac{1}{4} \xi_N^2 \cdot \frac{(1-z)^2}{1+z} + b^* z^2 \right)}{(1+z)^{1/2}}
$$
and thus
\begin{multline*}
\tilde{h}'(z) = \bigg(\frac{\frac{1}{2} \xi_N^2 \cdot \frac{1-z}{1+z} + \frac{1}{4} \xi_N^2 \cdot \frac{(1-z)^2}{(1+z)^2} + 2 b^* z}{(1+z)^{1/2}} - \frac{1/2}{(1+z)^{3/2}} \bigg) \\ \,\cdot\, \exp \left( - \tfrac{1}{4} \xi_N^2 \cdot \frac{(1-z)^2}{1+z} + b^* z^2 \right) \,.
\end{multline*}
Let
$$
Z \equiv Z_N := \big\{ z \in \mathbb{C} \,|\, z=re^{i\varphi}, 1-1/N \leq r \leq 1, |\varphi| \leq 1/\sqrt{N} \big\}
$$
and note that for $z \in Z$, we have the estimates 
\begin{align*}
\re \left( - \frac{(1-z)^2}{1+z} \right) \leq 4/N
\end{align*}
and thus
\begin{align*}
      \left| \tilde{h}'(z) \right| 
&\leq K \left\{ \Big( \xi_N^2 |1-z| + |b^*| + 1 \Big) \, \exp \left( \tfrac{1}{4} \xi_N^2 \cdot \re \bigg( - \frac{(1-z)^2}{1+z} \bigg) + |b^*| \right) \right\} \\
&\leq K \left\{ \Big( \xi_N^2 / \sqrt{N} + |b^*| + 1 \Big) \, \exp \left( \xi_N^2 / N + |b^*| \right) \right\} \,,
\end{align*}
where $K$ denotes some absolute constant.
Since for $z=Re^{it}$ with $|t| \leq 1/\sqrt{N}$,
the~line~segment between the points $z$~and~$1$
is obviously contained in the set $Z$, we~therefore~obtain
\begin{multline*}
\left| \tilde{h}(z) - \tilde{h}(1) \right|
\leq |z-1| \sup_{\alpha \in [0;1]} \left| \tilde{h}'((1-\alpha)z + \alpha) \right|
\leq |z-1| \sup_{\zeta \in Z} \left| \tilde{h}'(\zeta) \right|
\\
\leq K |z-1| \left( \xi_N^2/\sqrt{N} + |b^*| + 1 \right) \, \exp \Big( \xi_N^2 / N + |b^*| \Big)
\leq K(b^*,\xi^*) \, \sqrt{N} \, |z-1| \,,
\end{multline*}
where the last step uses the assumption that $\lim_{N \to \infty} \xi_N / \sqrt{N} = \xi$,
and $K(b^*,\xi^*)$ denotes some~constant depending only on $b^*$ and $\xi^* := (\xi_N)_{N \in \mathbb{N}}$.
Using (\ref{xi-est}), (\ref{est5}) and (\ref{theestimate}), it follows that
\begin{align*}
      |I_3|
&\leq \frac{1}{2\pi} \int_{-1/\sqrt{N}}^{+1/\sqrt{N}} \frac{\exp(|\eta|^2)}{|1-Re^{it}|^{5/2}} \cdot \left| \tilde{h}(Re^{it}) - \tilde{h}(1) \right| \frac{dt}{|Re^{it}|^N} \\
&\leq K(b^*,\xi^*,\eta) \, \sqrt{N} \int_{-1/\sqrt{N}}^{+1/\sqrt{N}} \frac{1}{|1-Re^{it}|^{3/2}} \ dt \\
&\leq K(b^*,\xi^*,\eta) \, N \left( 1 + \int_{1}^{\sqrt{N}\varepsilon} \frac{1}{u^{3/2}} \ du \right) \\
&\leq K(b^*,\xi^*,\eta) \, N
   =  o(N^{3/2}) \,,
\end{align*}
where $K(b^*,\xi^*,\eta)$ denotes some constant which depends only on $b^*$, $\xi^*$, and $\eta$ 
(and~which may change from line to line as usual).

For the integral $I_4$, we can finally use (\ref{xi-est}), (\ref{est5})
and (\ref{theestimate}) to obtain
\begin{align*}
      |I_4| 
&\leq \frac{1}{2\pi} \int_{1/\sqrt{N}}^{2\pi-1/\sqrt{N}} \frac{\exp(|\eta|^2 + |b^*|)}{|1-Re^{it}|^{5/2}} \ \frac{dt}{R^N} \\
&\leq K \exp(|\eta|^2 + |b^*|) \cdot \int_{1/\sqrt{N}}^{\pi} \frac{1}{|1-Re^{it}|^{5/2}} \ dt \\
&\leq K \exp(|\eta|^2 + |b^*|) \cdot N^{3/2} \int_{\sqrt{N}\varepsilon}^{\infty} \frac{1}{u^{5/2}} \ du \\
&\leq K \exp(|\eta|^2 + |b^*|) \cdot N^{3/4} = o(N^{3/2}) \,,
\end{align*}
where $K$ denotes some absolute constant which may change from step to step. 

This concludes the proof of Proposition \ref{proposition-B}.
\end{proof}

\bigskip

\end{document}